\documentclass{article}
\usepackage{amssymb}
\usepackage{amsmath}
\usepackage{amsfonts}

\newtheorem{theorem}{Theorem}

\newtheorem{definition}[theorem]{Definition}

\newtheorem{lemma}[theorem]{Lemma}

\newtheorem{remark}[theorem]{Remark}

\newenvironment{proof}[1][Proof]{\noindent\textbf{#1.} }{\ \rule{0.5em}{0.5em}}
\input{tcilatex}

\begin{document}

\title{Regular Variation and Smile Asymptotics}
\author{Shalom Benaim and Peter Friz \\
%EndAName
Statistical Laboratory, University of Cambridge}
\maketitle

\begin{abstract}
We consider risk-neutral returns and show how their tail asymptotics
translate \textit{directly} to asymptotics of the implied volatility smile,
thereby sharpening Roger Lee's celebrated moment formula. The theory of
regular variation provides the ideal mathematical framework to formulate and
prove such results. The practical value of our formulae comes from the vast
literature on tail asymptotics and our conditions are often seen to be true
by simple inspection of known results.
\end{abstract}

\section{\protect\bigskip Introduction}

Consider\ risk-neutral returns $X$ with cummulative distribution function $F$%
. We will impose a mild integrability condition on the right and left tail
denoted by (IR), (IL) respectively. We write $\bar{F}=1-F$ and, if it
exists, $f$ for the probability density function of $X$. $\,$The class of
regularly varying functions at $+\infty $ of index $\alpha $ is denoted by $%
R_{\alpha }$. The reader not familiar with the theory of regular variation
should think of positive multiples of $x^{\alpha }$ and harmless
perturbations such as $\left( 1+x^{\alpha /2}\right) ^{2},\,x^{\alpha }\log
x $ or $x^{\alpha }/\left( \log x\right) ^{3}$. By convention, functions in $%
R_{\alpha }$ are (eventually) positive.

The purpose of this paper is to connect the tail behaviour of $X$ to the
wing behaviour of the Black-Scholes implied volatility, sharpening Roger
Lee's celebrated moment formula \cite{Lee}. From a mathematical point of
view, the challenge is to relate the asymptotics of the distribution $F$ to
the asymptotics of a \textit{non-linear transform}, namely the Black-Scholes
implied volatility. From a financial point of view, we give further\
justification of volatility smile parametrizations seen in the industry and
obtain new insights into the wing behaviour of a variety of models.\newpage

The normalized price of a Black-Scholes call with log-strike $k$ is given by%
\begin{equation*}
c_{BS}\left( k,\sigma \right) =\Phi \left( d_{1}\right) -e^{k}\Phi \left(
d_{2}\right) \,
\end{equation*}%
with $d_{1,2}\left( k\right) =-k/\sigma \pm \sigma /2$. The implied
volatility is the (unique) value $V\left( k\right) $ so that 
\begin{equation*}
c_{BS}\left( k,V\left( k\right) \right) =\int_{k}^{\infty }\left(
e^{x}-e^{k}\right) dF\left( x\right) =:c\left( k\right) \text{.}
\end{equation*}%
We define the strictly decreasing function $\psi :\left[ 0,\infty \right]
\rightarrow \left[ 0,2\right] $ by 
\begin{equation*}
\psi \left[ x\right] =2-4\left[ \sqrt{x^{2}+x}-x\right] .
\end{equation*}

\begin{theorem}[Right-tail-wing formula]
\label{RTWF}Assume $\alpha >0$ and%
\begin{equation}
\exists \epsilon >0:\mathbb{E}[e^{(1+\epsilon )X}]<\infty .  \tag{IR}
\end{equation}%
Then%
\begin{equation*}
\text{(i)}\implies \text{(ii)}\implies \text{(iii)}\Longrightarrow \text{(iv)%
},
\end{equation*}%
where%
\begin{equation}
-\log f\left( k\right) \in R_{\alpha };  \tag{i}
\end{equation}%
\begin{equation}
-\log \bar{F}\left( k\right) \in R_{\alpha };  \tag{ii}
\end{equation}%
\begin{equation}
-\log c\left( k\right) \in R_{\alpha };  \tag{iii}
\end{equation}%
and\footnote{$g\left( k\right) \sim h\left( k\right) $ means $g\left(
k\right) /h\left( k\right) \rightarrow 1$as $k\rightarrow \infty .$}%
\begin{equation}
V(k)^{2}/k\sim \psi \left[ -\log c\left( k\right) /k\right] .\newline
\tag{iv}
\end{equation}%
\newline
If (ii) holds then $-\log c\left( k\right) \sim -k-\log \bar{F}$ and 
\begin{equation}
V(k)^{2}/k\sim \psi \left[ -1-\log \bar{F}\left( k\right) /k\right] , 
\tag{iv'}
\end{equation}%
if (i) holds, then $-\log f\sim -\log \bar{F}$ and%
\begin{equation}
V(k)^{2}/k\sim \psi \left[ -1-\log f\left( k\right) /k\right] . 
\tag{iv\textquotedblright }
\end{equation}%
Finally, if either $-\log f\left( k\right) /k$ or $-\log \bar{F}\left(
k\right) /k$ or $-\log c\left( k\right) /k$ goes to infinity as $%
k\rightarrow \infty $ then $V^{2}\left( k\right) $ behaves sublinearly. More
precisely,%
\begin{equation}
V(k)^{2}/k\sim \frac{1}{-2\log f\left( k\right) /k}\text{or }\frac{1}{-2\log 
\bar{F}\left( k\right) /k}\text{ or }\frac{1}{-2\log c\left( k\right) /k}. 
\tag{v}
\end{equation}
\end{theorem}

We emphasize that (iv),(iv'),(iv\textquotedblright ) contain the \textit{%
full asymptotics} of the implied volatility smile. For instance, we can see
when $\lim \sup V\left( k\right) ^{2}/k$ in Lee's moment formula is a
genuine limit: $V\left( k\right) ^{2}/k$ converges if and only if $-\log 
\bar{F}\left( k\right) /k$ converges to some limit $\theta $. Note that our
condition (IR) forces $\theta >1$. Note also that in this case $-\log \bar{F}%
\in R_{1}$ so that condition (ii) is automatically satisfied.

In models without moment explosion (Black-Scholes, Merton's jump diffusion
model, FMLS with $\beta =-1$...) the moment formula indicates sublinear
behaviour of the implied variance, $\lim \sup V\left( k\right) ^{2}/k=\lim
V\left( k\right) ^{2}/k=0$, but yields no further information. In contrast,
theorem \ref{RTWF} gives the \textit{precise sublinear asymptotics. }For
instance, in the Black-Scholes model $\log f_{BS}\left( k\right) \sim
-k^{2}/\left( 2\sigma ^{2}\right) $ and (v) implies $V\left( k\right)
^{2}\sim \sigma ^{2}$ in trivial agreement with the Black-Scholes flat
volatility smile.

There is a similar result which, as $k\rightarrow \infty $, links $f\left(
-k\right) ,F\left( -k\right) $, normalized out-of-the-money put prices $%
p\left( -k\right) $ and the implied volatility in the left wing.

\begin{theorem}[Left-tail-wing formula]
\label{LTWF}\label{RTWF copy(1)}Assume $\alpha >0$ and%
\begin{equation}
\exists \epsilon >0:\mathbb{E}[e^{-\epsilon X}]<\infty .  \tag{IL}
\end{equation}%
Then%
\begin{equation*}
\text{(i)}\implies \text{(ii)}\implies \text{(iii)}\Longrightarrow \text{(iv)%
},
\end{equation*}%
where%
\begin{equation}
-\log f\left( -k\right) \in R_{\alpha };  \tag{i}
\end{equation}%
\begin{equation}
-\log F\left( -k\right) \in R_{\alpha };  \tag{ii}
\end{equation}%
\begin{equation}
-\log p\left( -k\right) \in R_{\alpha };  \tag{iii}
\end{equation}%
and%
\begin{equation}
V(-k)^{2}/k\sim \psi \left[ -1-\log p\left( -k\right) /k\right] .\newline
\tag{iv}
\end{equation}%
\newline
If (ii) holds then $-\log p\left( -k\right) \sim k-\log F\left( -k\right) $
and 
\begin{equation}
V(-k)^{2}/k\sim \psi \left[ -\log F\left( -k\right) /k\right] ,  \tag{iv'}
\end{equation}%
if (i) holds, then $-\log f\left( -k\right) \sim -\log F\left( -k\right) $
and%
\begin{equation}
V(k)^{2}/k\sim \psi \left[ -\log f\left( -k\right) /k\right] . 
\tag{iv\textquotedblright }
\end{equation}%
Finally, if either $-\log f\left( -k\right) /k$ or $-\log F\left( -k\right)
/k$ or $-\log p\left( -k\right) /k$ goes to infinity as $k\rightarrow \infty 
$ then $V^{2}\left( -k\right) $ behaves sublinearly. More precisely,%
\begin{equation}
V(-k)^{2}/k\sim \frac{1}{-2\log f\left( -k\right) /k}\text{or }\frac{1}{%
-2\log F\left( -k\right) /k}\text{ or }\frac{1}{-2\log p\left( -k\right) /k}.
\tag{v}
\end{equation}
\end{theorem}

In conclusion, under mild integrability and regular variation conditions,
tail asymptotics translate \textit{directly} to asymptotics of the implied
volatility smile. The practical value of formulae (iv'),
(iv\textquotedblright ) comes from the vast literature on tail asymptotics
and our conditions are often seen to be true by simple inspection of known
results. The authors would like to thank Jim Gatheral, Roger Lee and Chris
Rogers for related discussions. Financial support from the Cambridge
Endowment for Research in Finance (CERF)\ is gratefully acknowledged.

\section{Elements of Regular Variation Theory}

\begin{definition}
\label{DefRV}A positive real-valued measurable function $f$ is regularly
varying with index $\alpha $, in symbols $g\in R_{\alpha }$ if 
\begin{equation*}
\lim_{x\rightarrow \infty }\frac{g(\lambda x)}{g(x)}=\lambda ^{\alpha }.
\end{equation*}%
Functions in $R_{0}$ are called slowly varying.
\end{definition}

The definition extends immediately to functions which are positive for $x$
large enough. Recall that $g\sim h$ means $\lim_{x\rightarrow \infty
}g(x)/h\left( x\right) =1$. Trivially, $g\in R_{\alpha }$ and $h\sim g$
implies $h\in R_{\alpha }$. The following result can be found in the fine
monograph \cite[Thm 4.12.10, p255]{BGT}.

\begin{theorem}[Bingham's Lemma]
\label{RegVarThm} Let $g\in R_{\alpha }$ with $\alpha >0$ such that that $%
e^{-g}$ is locally integrable at $+\infty $. Then%
\begin{equation*}
-\log \int_{x}^{\infty }e^{-g\left( y\right) }dy\sim g\left( x\right) .
\end{equation*}
\end{theorem}

\section{Right-Wing Smile Asympotics}

\begin{proof}[Proof of Theorem \protect\ref{RTWF}]
We first remark in presence of condition (IR), under either assumption (ii)
or (iii) one must actually have $\alpha \geq 1$. Indeed,%
\begin{equation*}
\bar{F}\left( k\right) =\mathbb{P}[X>k]\leq e^{-\left( 1+\epsilon \right) k}%
\mathbb{E}[e^{(1+\epsilon )X}]\implies -\log \bar{F}\left( k\right) \geq
\left( 1+\epsilon \right) k.
\end{equation*}%
Assuming $-\log \bar{F}\in R_{\alpha }$ for $\alpha <1$ would entail $-\log 
\bar{F}\left( k\right) \leq k^{\alpha ^{\prime }}$ for $\alpha ^{\prime }\in
\left( \alpha ,1\right) $ and $k$ large enough which contradicts the lower
bound just established. Similarly, 
\begin{equation*}
c\left( k\right) =\mathbb{E}[(e^{X}-e^{k})^{+}]\leq \mathbb{E}%
[e^{X};X>k]\leq \mathbb{E}[e^{X}e^{\epsilon (X-k)}]=e^{-\epsilon k}\mathbb{E}%
[e^{(1+\epsilon )X}]
\end{equation*}%
so that $-\log c\left( k\right) \geq \epsilon k$ and $-\log c\in R_{\alpha }$
can only happen with $\alpha \geq 1$.\newline
(i)$\implies $(ii): Use Binghan's lemma with the (eventually) positive
function $g\left( k\right) =-\log f\left( k\right) \in R_{\alpha }.$ (ii)$%
\implies $(iii): Lemma \ref{CallRegVar}. (iii)$\implies $(iv): Lemma \ref%
{Impvol/call}.The final statement on sublinear asymptotics follows from $%
\psi \left[ x\right] \sim 1/\left( 2x\right) $ as $x\rightarrow \infty $.
Indeed,%
\begin{equation*}
\psi \left[ x\right] =2-4x\left( \sqrt{1+1/x}-1\right) =2-4x\left( \frac{1}{%
2x}-\frac{1}{8x^{2}}+O\left( x^{-3}\right) \right) \sim \frac{1}{2x}.
\end{equation*}
\end{proof}

\begin{lemma}
\label{CafterIBP}Assume (IR). Then the normalized call price with log-strike 
$k$ is given by 
\begin{equation}
c\left( k\right) =\int_{k}^{\infty }{e^{x}\bar{F}(x)dx}.  \label{CallFormula}
\end{equation}
\end{lemma}

\begin{proof}
Note that $\bar{F}=o\left( e^{-x}\right) $ since $\bar{F}\left( x\right) =%
\mathbb{P}[X>x]\leq e^{-\left( 1+\epsilon \right) x}\mathbb{E}%
[e^{(1+\epsilon )X}]$. Integration by parts gives%
\begin{eqnarray*}
c\left( k\right) &=&-\int_{k}^{\infty }\left( e^{x}-e^{k}\right) d\bar{F}%
\left( x\right) \\
&=&-\bar{F}\left( x\right) \left( e^{x}-e^{k}\right) |_{k}^{\infty
}+\int_{k}^{\infty }e^{x}\bar{F}\left( x\right) dx
\end{eqnarray*}%
and the boundary contribution disappears.
\end{proof}

\begin{lemma}
\label{CallRegVar}Assume (IR) and that $-\log {\bar{F}}\in {R}_{\alpha }$
for some $\alpha \geq 1$. Then $-\log c\in $ ${R}_{\alpha }$ and as $%
k\rightarrow \infty $,%
\begin{equation}
-\log c\left( k\right) \sim -k-\log \bar{F}\left( k\right) .\newline
\label{CallBehaviour}
\end{equation}
\end{lemma}

\begin{proof}
Obviously $k\mapsto k\in R_{1}$ and $-\log \bar{F}\in R_{\alpha },\alpha
\geq 1$. We want to apply theorem \ref{RegVarThm} with $k\mapsto $ $\varphi
\left( k\right) \equiv -\log \bar{F}\left( k\right) -k$ to obtain (\ref%
{CallBehaviour}) but need to be careful since in general the difference of
two regularly varying functions may not be regularly varying. From the last
proof we know that $-\log {\bar{F}(k)\geq }\left( 1+\epsilon \right) k$.
When $\alpha >1$ then $\varphi \left( k\right) $ is immediately seen to be
in $R_{\alpha }$. When $\alpha =1$ we can write $-\log \bar{F}\left(
k\right) =k\,L\left( k\right) $ with slowly varying $L\left( k\right) \geq
1+\epsilon $. To see that $\varphi \in R_{1}$ it suffices to write 
\begin{equation*}
\frac{\varphi \left( \lambda k\right) }{\varphi \left( k\right) }=\frac{%
L\left( \lambda k\right) -1}{L\left( k\right) -1}=1+\frac{L\left( \lambda
k\right) /L\left( k\right) -1}{1-1/L\left( k\right) }
\end{equation*}%
so that%
\begin{equation*}
\left\vert \frac{\varphi \left( \lambda k\right) }{\varphi \left( k\right) }%
-1\right\vert \leq \left\vert \frac{L\left( \lambda k\right) /L\left(
k\right) -1}{1-1/L\left( k\right) }\right\vert \leq \frac{1+\epsilon }{%
\epsilon }\left\vert L\left( \lambda k\right) /L\left( k\right) -1\right\vert
\end{equation*}%
and this tends to zero as $k\rightarrow +\infty $ since $L\in R_{0}$. In
either case, we can apply theorem \ref{RegVarThm} with $g=\varphi $ and
obtain 
\begin{equation*}
-\log c\left( k\right) \overset{(\ref{CallFormula})}{=}-\log
\int_{k}^{\infty }e^{x+\log \bar{F}\left( x\right) }dx\sim -\left( k+\log 
\bar{F}\left( k\right) \right) =\varphi \left( k\right) .
\end{equation*}
Note that $\varphi \left( \cdot \right) $ was seen to be in $R_{\alpha }$ so
that $-\log c\left( k\right) \sim \varphi \left( k\right) $ must also be in $%
R_{\alpha }$.
\end{proof}

\begin{lemma}
\label{Impvol/call}Assume (IR) and $-\log {c\in R}_{\alpha }$ for some $%
\alpha \geq 1$. Then, as $k\rightarrow \infty ,$ 
\begin{equation*}
\frac{\log c\left( k\right) }{k}=-\frac{k}{2V(k)^{2}}+\frac{1}{2}-\frac{%
V(k)^{2}}{8k}+O\left( \frac{\log k}{k}\right)
\end{equation*}%
and 
\begin{equation*}
\frac{V(k)^{2}}{k}\sim \psi \left( \frac{-\log {c\left( k\right) }}{k}%
\right) .
\end{equation*}
\end{lemma}

\begin{proof}
Appendix.
\end{proof}

\begin{remark}
In the preceding lemma the condition of regular variation can be replaced by
the weaker%
\begin{equation*}
\exists n\in \mathbb{N}:\lim \inf_{k\rightarrow \infty }V(k)k^{n}>1.
\end{equation*}
\end{remark}

\section{Left-Wing Smile Asympotics}

The proofs are similar and therefore omitted.

\begin{lemma}
Assume (IL). Then the normalized put price with log-strike $-k$ is given by%
\begin{equation*}
p\left( -k\right) =\int_{-\infty }^{-k}{e^{x}F(x)dx}=\int_{k}^{\infty
}e^{-x}F\left( -x\right) dx.
\end{equation*}
\end{lemma}

\begin{lemma}
Assume (IL) and that $-\log {F}\left( -k\right) $ is regularly varying. Then 
$-\log p\left( -k\right) $ is regularly varying as $k\rightarrow \infty $
and 
\begin{equation}
-\log p\left( -k\right) \sim k-\log {F(-k).}  \label{PutBehaviour}
\end{equation}
\end{lemma}

\begin{lemma}
Assume (IL) and that $-\log p\left( -k\right) $ is regularly varying. Then,
as $k\rightarrow \infty ,$ 
\begin{equation*}
\frac{V(-k)^{2}}{k}\sim \psi \left( -1+\frac{-\log p\left( -k\right) }{k}%
\right) .
\end{equation*}
\end{lemma}

\section{Examples}

In practice, one has $X=\log \left( S_{T}/F_{T}\right) $ where $S_{T}$
denotes the risk-neutral stock price at time $T$ and $F_{T}$ is the time-$T$
forward price. Also, $k=\log \left( K/F_{T}\right) $ and all quanitities $%
f,F,c,p,$\thinspace $V$ depend on time $T$ and we set $V\left( k\right)
=V\left( k,T\right) =:\sigma _{BS}\left( k,T\right) \sqrt{T}$, now calling $%
\sigma _{BS}\left( k,T\right) $ the implied volatility. It is worthwhile to
spell of part (iv) of Theorem \ref{RTWF},%
\begin{equation*}
\sigma _{BS}^{2}\left( k,T\right) T/k\sim \psi \left[ -\log c\left(
k,T\right) /k\right] .
\end{equation*}%
Noting that $\left( S_{T}/F_{T}:T\geq 0\right) $ is a martingale and using
convexity of the call payoff, it is easy to see, for $k$ fixed, $c\left(
k,T\right) =\mathbb{E}\left( e^{X}-e^{k}\right) ^{+}$ is non-decreasing in $%
T $ and so is $\psi \left[ -\log c\left( k,T\right) /k\right] $. In fact, J.
Gatheral and E. Reiner\footnote{%
Presentations at the Global Derivatives \& Risk Management Conference 2004.}
point out independently that the \textit{implied total variance} $\sigma
_{BS}^{2}\left( k,T\right) T$ \ is non-decreasing in $T$, a consequence of
monotonicity of the undiscounted Black-Scholes prices in $\sigma ^{2}T$.
Thus, our asymptotic results respect the term structure of the volatility
surface.\newline
In the examples below, we will focus mainly on applications of the \textit{%
right-tail-wing formula}, applications of the \textit{left-tail-wing formula}
being nearly identical.

\subsection{Sanity Check: Black-Scholes Model}

If $\sigma $ denotes the Black Scholes volatility, the returns have a normal
density with variance $\sigma ^{2}T.$ Obviously then,%
\begin{equation*}
\log f_{BS}\left( k\right) \sim -k^{2}/\left( 2\sigma ^{2}T\right)
\end{equation*}%
and Theorem \ref{RTWF} implies $\sigma _{BS}\left( k,T\right) ^{2}\sim
\sigma ^{2}$ as $k\rightarrow \infty $ in trivial agreement $V\equiv \sigma $%
.

\subsection{Barndorff-Nielsen's NIG Model}

Here $X=X_{T}\sim $ $NIG\left( \alpha ,\beta ,\mu T,\delta T\right) $. The
moment generating function is given by%
\begin{equation*}
M\left( z\right) =\exp \left[ T\left( \delta \left\{ \sqrt{\alpha ^{2}-\beta
^{2}}-\sqrt{\alpha ^{2}-\left( \beta +z\right) ^{2}}\right\} +\mu z\right) %
\right] .
\end{equation*}%
It is custom to write $\alpha =\sqrt{\beta ^{2}+\gamma ^{2}}$ with $\gamma
>0 $. From \cite{BN} and the references therein we have%
\begin{equation*}
f\left( k\right) \sim C\left\vert k\right\vert ^{-3/2}e^{-\sqrt{\beta
^{2}+\gamma ^{2}}\left\vert k\right\vert +\beta k}\text{ as }k\rightarrow
\pm \infty
\end{equation*}%
and we see that $-\log f$ is regularly varying (with index $1$). Moreover,%
\begin{equation*}
\log f\left( k\right) /k\rightarrow \left( -\sqrt{\beta ^{2}+\gamma ^{2}}%
+\beta \right) \text{ as }k\rightarrow +\infty \text{.}
\end{equation*}%
and from Theorem \ref{RTWF} 
\begin{eqnarray*}
\frac{\sigma _{BS}^{2}\left( k,T\right) T}{k} &\sim &\psi \left( -1-\log
f\left( k\right) /k\right) \\
&\sim &\psi \left( -1+\sqrt{\beta ^{2}+\gamma ^{2}}-\beta \right)
\end{eqnarray*}%
in agreement with Lee's moment formula with critical moment $1+p^{\ast }=$ $%
\sqrt{\beta ^{2}+\gamma ^{2}}-\beta -1=\alpha -\beta $, as can be seen
directly from the moment generating function, see \cite{Lee}.

\subsection{Carr-Wu's Finite Moment Logstable Model\label{CWexample}}

Here $X=X_{T}\sim L_{\alpha }\left( \mu T,\sigma T^{1/\alpha },-1\right) $
where the law $L_{\alpha }\left( \theta ,\sigma ,\beta \right) $ has
characteristic function%
\begin{equation*}
\mathbb{E}\left[ e^{iuX}\right] =e^{iu\theta -\left\vert u\right\vert
^{\alpha }\sigma ^{\alpha }\left( 1-i\beta \left( \text{sig }u\right) \tan 
\frac{\pi \alpha }{2}\right) .}
\end{equation*}%
with $\alpha \in (1,2],\theta \in \mathbb{R},\sigma \geq 0,\beta \in \lbrack
-1,1]$. From \cite[page 10, equation\ (6)]{CW} and the references therein,%
\begin{equation}
-\log \bar{F}\left( k\right) \sim k^{\frac{\alpha }{\alpha -1}}\times \left[
T\alpha \sigma ^{\alpha }\left\vert \sec \left( \pi \alpha /2\right)
\right\vert \right] ^{-1/\left( \alpha -1\right) }\,\text{as }k\rightarrow
\infty ,  \label{CarrWuTail}
\end{equation}%
and from Theorem \ref{RTWF} we see that%
\begin{equation*}
\sigma _{BS}^{2}\left( k,T\right) T\sim k^{1-\frac{1}{a-1}}\times \frac{1}{2}%
\left[ T\alpha \sigma ^{\alpha }\left\vert \sec \left( \pi \alpha /2\right)
\right\vert \right] ^{1/\left( \alpha -1\right) }\rightarrow 0\text{ as }%
k\rightarrow \infty \text{.}
\end{equation*}%
Note that in the limit $\alpha \uparrow 2$ the Black Scholes result is
recovered. We also note that the moment generating function of $L_{\alpha
}\left( \theta ,\sigma ,-1\right) \mathbb{\ }$exists for all positive $z$
and is given by 
\begin{equation*}
M\left( z\right) =\exp \left[ z\theta -\left( z\sigma \right) ^{\alpha }\sec
\left( \frac{\pi \alpha }{2}\right) \right] ;
\end{equation*}%
Kasahara's exponential Tauberian theorem \cite[p253]{BGT} then gives
immediately (\ref{CarrWuTail}).

\subsection{Merton's Jump Diffusion Model\label{MertonEx}}

As in the examples above, the return process $X_{\cdot }$ is L\'{e}vy with
triplet $(\mu ,\sigma ^{2},K)$ where $K$ is $\lambda $ (=intensity of jump)
times a Gaussian measure with mean $\alpha $ and standard deviation $\delta $
describing the distribution of jumps. We now demonstrate how to proceed
without explicit knowledge of the asymptotic tail. Set $\,X=X_{T}$ and note
that $\mathbb{E}\left[ \exp \left( zX\right) \right] <\infty $ for all $z$.
Optimizing over $z$ we get the tail estimate%
\begin{equation*}
\bar{F}\left( k\right) =\mathbb{P}\left[ X>k\right] \leq \inf_{z}e^{-zk}%
\mathbb{E}\left[ \exp \left( zX\right) \right] =e^{K\left( z^{\ast }\right)
-z^{\ast }k}
\end{equation*}%
where $K\left( z\right) =\log \mathbb{E}\left[ \exp \left( zX\right) \right] 
$ is the logarithmic mgf and $z^{\ast }=z^{\ast }\left( k\right) $ is
determined from\footnote{%
Readers familiar with large deviation theory recognize the Fenchel-Legendre
transform.}%
\begin{equation*}
K^{\prime }\left( z^{\ast }\right) =k\text{.}
\end{equation*}%
For the Merton model,%
\begin{equation*}
K\left( z\right) =T\left\{ z\mu +\frac{1}{2}z^{2}\sigma ^{2}+\lambda \left(
e^{z\alpha +z^{2}\delta ^{2}/2}-1\right) \right\}
\end{equation*}%
from which for $\delta >0$ it is easy to see that%
\begin{equation*}
z^{\ast }=z^{\ast }\left( k\right) \sim \frac{\sqrt{2\log k}}{\delta }
\end{equation*}%
so that $K\left( z^{\ast }\right) -z^{\ast }k\sim -z^{\ast }k$ and%
\begin{equation*}
\log \bar{F}\left( k\right) \lesssim -z^{\ast }k=-\frac{k}{\delta }\sqrt{%
2\log k}.
\end{equation*}%
From the saddle point results of \cite{EJML} or the L\'{e}vy tail estimates
from \cite{AB2} specialized to this example, 
\begin{equation*}
\log \bar{F}\left( k\right) \sim -\frac{k}{\delta }\sqrt{2\log k},
\end{equation*}%
and Theorem \ref{RTWF} implies%
\begin{equation*}
\sigma _{BS}^{2}\left( k,T\right) T\sim \delta \times \frac{k}{2\sqrt{2\log k%
}}.
\end{equation*}%
Note that this is independent of the mean jump size $\alpha $ provided $%
\delta >0$. If $\delta =0$ and $\alpha >0$ a similar argument shows that $%
z^{\ast }=z^{\ast }\left( k\right) \sim \log k/\alpha $ and%
\begin{equation}
\sigma _{BS}^{2}\left( k,T\right) T\sim \alpha \times \frac{k}{2\log k}.
\label{MertonAsympDeltaZero}
\end{equation}

\begin{remark}[J. Gatheral]
A Poisson process with intensity $\lambda $ has $n$ jumps with probability $%
e^{-\lambda }\lambda ^{n}/n!=e^{g\left( n\right) }$ with $g\left( n\right)
\sim -n\log n$ by Stirling's formula. The time $T=1$ Black-Scholes value of
a digital is $\Phi \left( d_{2}\right) $ and $\log \Phi \left( d_{2}\right)
\sim -k^{2}/2\sigma ^{2}$. Identifying $k\sim n\alpha $ leads to%
\begin{equation*}
-k^{2}/2\sigma ^{2}\sim -\left( k/\alpha \right) \log k\implies \sigma
^{2}\sim \alpha \times \frac{k}{2\log k}
\end{equation*}%
in agreement with (\ref{MertonAsympDeltaZero}).
\end{remark}

\section{Further Examples and Discussion}

Our tail wing formulae apply as soon as one has some (in fact: crude)
asymptotic knowledge of the tail behaviour of the returns. As remarked in
the introduction, for many models fine tail asymptotics are available in the
literature and the smile asymptotics follow. In particular, the estimates of
Albin-Bengtsson \cite{AB1, AB2} allow to use our results for the vast
majority of exponential L\'{e}vy models \cite{S}: NIG appears as a special
case of generalized hyperbolic processes, Meixner processes as special cases
of GZ\ processes. CGMY and Variance Gamma model are also covered. Tail
estimates for stochastic volatility models appear in the literature,
although one has to be careful with results obtained in short time regimes
since the limits $T\rightarrow 0$ and $k\rightarrow \infty $ in general do
not interchange\footnote{%
The wrong wing behaviour in Hagan's SABR\ formula, based on short time
asymptotics, is a warning example.
\par
{}}. We remark that condition (IR) rules out models for which \textit{every} 
$p^{th}$-moment of the underlying, $p>1$, explodes and a similar remark
applies to condition (IL). Unfortunately, there are stochastic volatility
models with this degenerate behaviour \cite{AP} and our (current) results do
not apply.\newline
When a moment generating function is known, one can often use an exponential
Tauberian theorem to obtain log-tail estimates, as demonstrated in the
Example \ref{CWexample}. Also, the Fenchel-Legendre transform gives quick
upper bounds which often turn out to be sharp, compare Example \ref{MertonEx}%
. We finally remark that some of the implications in Theorem \ref{RTWF} and
Theorem \ref{LTWF} can be reversed under Tauberian conditions but we shall
not pursue this here.

\section{Appendix}

\begin{proof}[Proof of Lemma \protect\ref{Impvol/call}]
The following bounds for the normal distribution function $\Phi $ are
well-known and can be obtained by integration by parts (or other methods),%
\begin{equation*}
\frac{e^{-x^{2}/2}}{\sqrt{2\pi }x}(1-\frac{1}{x^{2}})\leq \Phi (-x)\leq 
\frac{e^{-x^{2}/2}}{\sqrt{2\pi }x},\,\,\,x>0\text{.}
\end{equation*}%
Using assumption (IR), 
\begin{equation*}
E[(e^{X}-e^{k})_{+}]\leq E[e^{X};X>k]\leq E[e^{X}e^{\epsilon
(X-k)}]=e^{-\epsilon k}E[e^{(1+\epsilon )X}],
\end{equation*}%
we see that $c\left( k\right) $ is exponentially small in $k$, and this
implies that 
\begin{equation}
a:=\limsup_{k\rightarrow \infty }\frac{V(k)^{2}}{k}<2  \label{limsupLT2}
\end{equation}%
In particular, both%
\begin{equation*}
d_{1}(k)=-\frac{k}{V(k)}+\frac{V(k)}{2},\,\,\,d_{2}(k)=-\frac{k}{V(k)}-\frac{%
V(k)}{2}
\end{equation*}%
will then be negative for $k$ large enough so that we will be able to use
the bounds on $\Phi $. From the Black-Scholes formula for a normalized call
with log-strike $k$,%
\begin{equation*}
c\left( k\right) =\Phi (d_{1})-e^{k}\Phi (d_{2}),
\end{equation*}%
we obtain%
\begin{equation*}
e^{-d_{1}^{2}/2}\left( -\frac{1}{d_{1}}\right) \left( 1-\frac{1}{d_{1}^{2}}%
\right) -e^{k}e^{-d_{2}^{2}/2}\left( -\frac{1}{d_{2}}\right) \leq \sqrt{2\pi 
}c\left( k\right)
\end{equation*}%
and%
\begin{equation*}
\sqrt{2\pi }c\left( k\right) \leq e^{-d_{1}^{2}/2}\left( -\frac{1}{d_{1}}%
\right) -e^{k}e^{-d_{2}^{2}/2}\left( -\frac{1}{d_{2}}\right) \left( 1-\frac{1%
}{d_{2}^{2}}\right) .
\end{equation*}%
Because $d_{1}^{2}/2=d_{2}^{2}/2-k$, this simplifies to%
\begin{equation*}
e^{-d_{1}^{2}/2}\left[ -\frac{1}{d_{1}}\left( 1-\frac{1}{d_{1}^{2}}\right) +%
\frac{1}{d_{2}}\right] \leq \sqrt{2\pi }c\left( k\right) \leq
e^{-d_{1}^{2}/2}\left[ -\frac{1}{d_{1}}+\frac{1}{d_{2}}\left( 1-\frac{1}{%
d_{2}^{2}}\right) \right] .
\end{equation*}%
We now define $\epsilon _{1}=\epsilon _{1}\left( k\right) $ by%
\begin{equation*}
\log c\left( k\right) =-\frac{d_{1}^{2}}{2}+\epsilon _{1}(k)
\end{equation*}%
and will show $\left\vert \epsilon _{1}\left( k\right) \right\vert =O\left(
\log k\right) $. To start, we note the bounds%
\begin{equation}
-\frac{1}{d_{1}}\left( 1-\frac{1}{d_{1}^{2}}\right) +\frac{1}{d_{2}}\leq 
\sqrt{2\pi }e^{\epsilon _{1}(k)}\leq -\frac{1}{d_{1}}+\frac{1}{d_{2}}\left(
1-\frac{1}{d_{2}^{2}}\right) .  \label{eps1bounds}
\end{equation}%
From (\ref{limsupLT2}) we know that $V(k)\leq \sqrt{2}k$, hence%
\begin{equation*}
-d_{2}=\frac{k}{V(k)}+\frac{V(k)}{2}\geq \sqrt{k/2}\rightarrow \infty \text{
as }k\rightarrow \infty .
\end{equation*}%
Therefore, for $k$ large enough,%
\begin{equation*}
\sqrt{2\pi }e^{\epsilon _{1}(k)}\leq -\frac{1}{d_{1}}.
\end{equation*}%
In fact, (\ref{limsupLT2}) gives $a^{\prime }\in \left( 0,2\right) $ s.t. $%
V\left( k\right) ^{2}<a^{\prime }k$ for all $k$ large enough. From the
definition of $d_{1}$ we then see that%
\begin{equation*}
a^{\prime \prime \prime }\sqrt{k}:=\sqrt{k}\left( \frac{1}{\sqrt{a^{\prime }}%
}-\frac{\sqrt{a^{\prime }}}{2}\right) <-d_{1}\left( k\right)
\end{equation*}%
with $a^{\prime \prime \prime }>0$ because $a^{\prime }\in (0,2)$. Hence%
\begin{equation*}
\sqrt{2\pi }e^{\epsilon _{1}(k)}\leq \frac{1}{a^{\prime \prime \prime }\sqrt{%
k}}
\end{equation*}%
and we see that $\epsilon _{1}\left( k\right) \rightarrow -\infty $ as $%
k\rightarrow \infty $ and we only need a bound on $-\epsilon _{1}\left(
k\right) $. We start by showing

\begin{equation}
\exists n\in \mathbb{N}:\lim \inf_{k\rightarrow \infty }V(k)k^{n}>1.
\label{liminf}
\end{equation}%
To see this note that for $k$ large enough $\epsilon _{1}<0$ so that 
\begin{equation*}
-\log c\left( k\right) =\frac{d_{1}^{2}}{2}-\epsilon _{1}(k)\geq \frac{%
d_{1}^{2}}{2}
\end{equation*}%
By assumption $-\log c\in R_{\alpha }$ for $\alpha \geq 1$ so there exists $%
L\in R_{0}$ so that $-\log c\left( k\right) =k^{\alpha }L\left( k\right) $
and 
\begin{equation*}
k^{\alpha }L\left( k\right) \geq \frac{d_{1}^{2}}{2}=\frac{1}{2}\left( -%
\frac{k}{V(k)}+\frac{V(k)}{2}\right) ^{2}\geq \frac{k^{2}}{2V\left( k\right)
^{2}}-\frac{k}{2}
\end{equation*}%
so that%
\begin{equation*}
\frac{k^{2}}{V\left( k\right) ^{2}}\leq 2k^{\alpha }L\left( k\right)
+k\equiv k^{\alpha }\tilde{L}\left( k\right)
\end{equation*}%
where $\tilde{L}\in R_{0}$ and $V\left( k\right) ^{2}k^{\alpha -2}=1/\tilde{L%
}\left( k\right) $ and (\ref{liminf}) holds with any integer $n>\left(
\alpha -2\right) ^{+}.$

We return to establish a bound on $-\epsilon _{1}\left( k\right) $. We
already have lower and upper bounds on the implied volatility, namely
(replace $n$ by $n+1$ if needed)%
\begin{equation*}
k^{-2n}<V\left( k\right) ^{2}<a^{\prime }k\text{ with }a^{\prime }<2\text{
and }k\text{ large enough,}
\end{equation*}%
which we can use to bound $d_{1,2}$, namely%
\begin{eqnarray}
-d_{1} &<&k^{1+n}-\frac{1}{2k^{n}}\leq k^{1+n}  \label{UpperboundOnMinusD12}
\\
-d_{2} &<&k^{1+n}+\frac{\sqrt{2k}}{2k}\leq 2k^{1+n}  \notag
\end{eqnarray}%
(at least for $k\geq 1\,$). In order to derive an upper bound for $-\epsilon
_{1}$ recall that a lower bound on $\sqrt{2\pi }e^{\epsilon _{1}}$ was given
in (\ref{eps1bounds}) by 
\begin{eqnarray*}
\epsilon _{2}(k) &\equiv &-\frac{1}{d_{1}}\left( 1-\frac{1}{d_{1}^{2}}%
\right) +\frac{1}{d_{2}} \\
&=&\frac{d_{1}-d_{2}}{d_{1}d_{2}}+\frac{1}{d_{1}^{3}}=\frac{V}{d_{1}d_{2}}+%
\frac{1}{d_{1}^{3}} \\
&=&\frac{V^{2}d_{1}^{2}+Vd_{2}}{Vd_{1}^{3}d_{2}} \\
&=&\frac{\left( -k+V^{2}/2\right) ^{2}+\left( -k-V^{2}/2\right) }{%
Vd_{1}^{3}d_{2}}=:\left( \ast \right)
\end{eqnarray*}%
We already noted the existence of $\alpha ^{\prime }\in \left( 0,2\right) $
such that $x:=V^{2}/2<\beta k<k$ with $\beta =\alpha ^{\prime }/2<1$ for $k$
large enough. Noting that, for $k$ fixed, the function%
\begin{equation*}
\lbrack 0,\beta k]\ni x\mapsto \left( -k+x\right) ^{2}+\left( -k-x\right)
\end{equation*}%
is strictly decreasing on $[0,\beta k]$ we can get a lower bound on $%
\epsilon _{2}$ as follows,%
\begin{eqnarray*}
\left( \ast \right) &\geq &\frac{\left( -k+\beta k\right) ^{2}+\left(
-k-\beta k\right) }{d_{1}^{3}d_{2}V} \\
&=&\frac{k^{2}\left( 1-\beta \right) ^{2}-k\left( 1+\beta \right) }{%
d_{1}^{3}d_{2}V} \\
&\geq &\frac{\frac{1}{2}k^{2}\left( 1-\beta \right) ^{2}}{d_{1}^{3}d_{2}%
\sqrt{\beta k}}\text{ \ (for }k\text{ large enough)} \\
&\geq &\frac{\frac{1}{2}k^{2}\left( 1-\beta \right) ^{2}}{\left(
-d_{1}\right) ^{3}\left( -d_{2}\right) \sqrt{\beta k}} \\
&=&\frac{\frac{1}{2}k^{2}\left( 1-\beta \right) ^{2}}{\left( k^{1+n}\right)
^{3}\left( 2k^{1+n}\right) \sqrt{\beta k}}\text{ \ (use (\ref%
{UpperboundOnMinusD12}))} \\
&=&\gamma k^{2-4\left( 1+n\right) -1/2}\text{ with }\gamma >0.
\end{eqnarray*}

Therefore,%
\begin{equation*}
\sqrt{2\pi }e^{\epsilon _{1}(k)}>\epsilon _{2}(k)>\gamma k^{-5/2-4n}>0
\end{equation*}

and after taking logarithms%
\begin{equation*}
-\epsilon _{1}(k)<\left( 4n+5/2\right) \log k+C
\end{equation*}%
for some constant $C=C\left( \gamma \right) $. In particular,%
\begin{equation*}
\left\vert \epsilon _{1}(k)\right\vert =O(\log k).
\end{equation*}%
Recall that $\epsilon _{1}$ was defined by%
\begin{equation*}
\log c\left( k\right) =-\frac{d_{1}^{2}}{2}+\epsilon _{1}(k)=\frac{k}{2}-%
\frac{k^{2}}{2V^{2}}-\frac{V^{2}}{8}+\epsilon _{1}(k)
\end{equation*}%
and after dividing by $k$ we have,%
\begin{equation*}
\left\vert \frac{\log c\left( k\right) }{k}-\left( \frac{1}{2}-\frac{k}{%
2V^{2}}-\frac{V^{2}}{8k}\right) \right\vert =O\left( \log k/k\right)
\end{equation*}%
which tends to zero as $k\rightarrow \infty $ so that%
\begin{equation*}
\frac{\log c\left( k\right) }{k}\sim \frac{1}{2}-\frac{k}{2V^{2}}-\frac{V^{2}%
}{8k}.
\end{equation*}%
We can rearrange the above equations as a quadratic equation in $V^{2}/k$,%
\begin{equation*}
-\frac{\log c\left( k\right) }{k}+\frac{\epsilon _{1}(k)}{k}=\frac{1}{%
2V^{2}/k}+\frac{V^{2}/k}{8}-\frac{1}{2}.
\end{equation*}%
One solution is given by%
\begin{equation*}
\frac{V^{2}}{k}=\psi \left( -\frac{\log c\left( k\right) }{k}+\frac{\epsilon
_{1}(k)}{k}\right) \sim \psi \left( -\frac{\log c\left( k\right) }{k}\right)
,
\end{equation*}%
while the other is rejected as it violates (\ref{limsupLT2}).
\end{proof}

\end{document}